\documentstyle[12pt]{article}
\begin{document}

\bibliographystyle{tom}

\newfont{\gothic}{eufm10 scaled\magstep1}
\newfont{\gothics}{eufm7 scaled\magstep1}

\newcommand{\gotA}{\mbox{\gothic A}}
\newcommand{\gotB}{\mbox{\gothic B}}
\newcommand{\gotC}{\mbox{\gothic C}}
\newcommand{\gotF}{\mbox{\gothic F}}
\newcommand{\gotI}{\mbox{\gothic I}}
\newcommand{\gotK}{\mbox{\gothic K}}
\newcommand{\gotM}{\mbox{\gothic M}}
\newcommand{\gotU}{\mbox{\gothic U}}
\newcommand{\sgotA}{\mbox{\gothics A}}
\newcommand{\sgotB}{\mbox{\gothics B}}
\newcommand{\sgotI}{\mbox{\gothics I}}
\newcommand{\sgotK}{\mbox{\gothics K}}
\newcommand{\sgotP}{\mbox{\gothics P}}

\newcommand{\Ii}{{\bf I}}
\newcommand{\Ni}{{\bf N}}
\newcommand{\Ri}{{\bf R}}
\newcommand{\Ci}{{\bf C}}
\newcommand{\Zi}{{\bf Z}}
\newcommand{\Ti}{{\bf T}}

\newcommand{\ch}{{\cal H}}

\newcommand{\proof}{\mbox{\bf Proof} \hspace{5pt}}
\newcommand{\remark}{\mbox{\bf Remark} \hspace{5pt}}

\newcommand{\ad}{{\mathop{\rm ad}}}
\newcommand{\Ad}{\mathop{\rm Ad}}
\newcommand{\RRe}{\mathop{\rm Re}}
\newcommand{\IIm}{\mathop{\rm Im}}
\newcommand{\Tr}{{\mathop{\rm Tr}}}
\newcommand{\supp}{\mathop{\rm supp}}

\newcommand{\disjoint}{\mathrel{\raisebox{-2.5pt}{$\circ
$\makebox[0pt]{\raisebox{1pt}{\hspace{-6pt}\hbox{\vrule height8pt width.4pt
depth-4pt}}}}}}

\renewcommand{\Box}{\raisebox{0.6ex}{\framebox[0.6em]{\rule{0em}{0.6ex}}}}

\newenvironment{remarkn}{\begin{rem} \rm}{\end{rem}}
\newenvironment{exam}{\begin{voorb} \rm}{\end{voorb}}

\newtheorem{lemma}{Lemma}[section]
\newtheorem{thm}[lemma]{Theorem}
\newtheorem{cor}[lemma]{Corollary}
\newtheorem{voorb}[lemma]{Example}
\newtheorem{rem}[lemma]{Remark}
\newtheorem{prop}[lemma]{Proposition}

\thispagestyle{empty}
\begin{center}

\vspace*{15mm}

{\Large{\bf Abundance of invariant and almost invariant}}  \\[2mm]
{\Large{\bf  pure states of $C^*$-dynamical systems}}  \\[4mm]

{\large{Ola Bratteli}}\\[1mm]

Department of Mathematics   \\
University of Oslo \\
Blindern, P.O.Box 1053, \\
N-0316 Oslo, Norway\\[2mm]

{\large{Akitaka Kishimoto}}\\[1mm]

Department of Mathematics \\
Hokkaido University \\
Sapporo 060, Japan\\[2mm]

and \\[2mm]

{\large{Derek W. Robinson}}\\[1mm]
Centre for Mathematics and its Applications\\
Australian National University\\
Canberra, ACT 0200,
Australia\\[10mm]

{\large{\bf Abstract}}
\begin{list}{}{\leftmargin=1.8cm \rightmargin=1.8cm}
\item
We show that invariant states of $C^*$-dynamical systems can
be approximated in the weak$^*$-topology by
invariant pure states, or  almost invariant pure states,
under various circumstances.

\end{list}

\end{center}

\newpage
\setcounter{page}{1}

\section{Introduction}

The states $E_{\sgotA}$ of a unital $C^*$-algebra $\gotA$
are the weak$^*$-closure of the convex envelope of the pure
states, i.e., the extremal states, and if $\alpha$ is a
$^*$-automorphism of $\gotA$ then the $\alpha$-invariant
states $E^\alpha_{\sgotA}$ are the weak$^*$-closure of the
convex envelope of the extremal $\alpha$-invariant states,
i.e., the extremal points of $E^\alpha_{\sgotA}$.
These statements are consequences of convexity and
compactness arguments.
In particular these arguments, and specifically the
Krein--Milman theorem, ensure the existence of an
abundant set of extremal and extremal $\alpha$-invariant states.
This type of reasoning does not, however, establish the
existence of  $\alpha$-invariant pure states and this is
not surprising as there are  examples for which no
invariant pure states exist.
Thus it is of interest to analyze the structural
properties of $(\gotA, \alpha)$ which ensure the
abundance of  invariant pure states.

The state space of a $C^*$-algebra often has much more
striking geometric properties.
If $\gotA$ is a UHF-algebra, or more generally if
$\gotA$ is anti-liminal, then the pure states are
weak$^*$-dense in $E_{\sgotA}$ (see \cite{Gli} or
 \cite{BR1}, Example 4.1.31).
It is not necessary to take convex envelopes.
Fannes, Nachtergaele and Werner \cite{FNW} have shown that
a similar situation can occur for invariant states.
If $\gotA=\otimes^\infty_{-\infty}M_n={M_n}^{\otimes\infty}$
is the UHF-algebra formed as the infinite tensor product of the
algebra $M_n$ of all $n\times n$-matrices and $\alpha$ is the shift
automorphism then it follows from \cite{FNW} that the
$\alpha$-invariant pure states are dense in $E^\alpha_{\sgotA}$.
The proof of this result relies on the construction of a
special class of states, the finitely-correlated states,
and appears limited to this particular example.
In particular it appears difficult to extend the result
to the shift action of $\Zi^d$ on $\otimes_{\Zi^d}M_n$
when $d>1$ and the question of density of the invariant pure states
in this situation appears to remain open.
This contrasts to the well known result that the extremal invariant
states are weak$^*$-dense in the invariant states for this example
(see \cite{BR1}, Example 4.3.26).
In this paper we will demonstrate that this rather
surprising density property of the invariant pure states,
and related properties, are true for some  other classes of
$C^*$-dynamical system.

First, in Section~2, we establish a slightly different but related
density result.
Under general assumptions on the $C^*$-dynamical system
$(\gotA, G, \alpha)$ we demonstrate density of the
`almost' invariant pure states.
The class of systems covered by our assumptions contains the
higher dimensional spin systems mentioned in the previous paragraph.
It also includes systems where the set of  invariant pure states is
not weak$^*$-dense in the set of invariant states and even systems
for which the set of invariant pure states is weak$^*$-closed or
even empty.
In Section~3 we then examine systems for which the group $G$ is
compact and demonstrate density of the  invariant pure states for a
broad class of algebras whenever the action of the group is a
quasi-product action \cite{BKR}, \cite{BEEK}, \cite{BEK}.
These actions encompass product type actions of $G$ on
UHF-algebras.
Finally, in Section~4, we extend the
Fannes--Nachtergaele--Werner result to the shift automorphism acting
on the infinite tensor product of a prime, unital $AF$-algebra.


\section{Approximation of invariant states by almost invariant pure
states  for discrete amenable group actions}

If $\gotA$ is a $C^*$-algebra let $\gotA_{\sgotP}$ denote the
set of positive elements $e\in\gotA$ such that $0\leq e\leq 1$ and
there is a positive non-zero element $p\in\gotA$ with $p\,e=p$.
By spectral modification it is easy to show that the element $p$
can be taken from $\gotA_{\sgotP}$ if it exists and any
positive $e\in\gotA$ with $\|e\|=1$ can be approximated in norm by
elements from $\gotA_{\sgotP}$.
These elements will play the role of projections in $\gotA$ if
$\gotA$ is not of real rank zero.

Recall from \cite{Ell}, \cite{Kis1} that an automorphism $\alpha$ of
a $C^*$-algebra $\gotA$ is said to be properly outer if for any
$y\in\gotA_{\sgotP}$   and any
$\varepsilon>0$ there exists an $x\in\gotA_{\sgotP}$ with
$x\,y=x$  such that
\[
\|x\,\alpha(x)\|<\varepsilon\;\;\;.
\]
A long list of equivalent conditions for proper outerness, when
$\gotA$ is separable, is given in Theorem 6.6 in \cite{OlP}.
A particularly useful characterization is given in Theorem~2.1 of
\cite{Kis2}: $\alpha$ is properly outer if for any
$y\in\gotA_{\sgotP}$, any $a\in\gotA\bigcup\{\Ii\}$ and any
$\varepsilon>0$ there exists an $x\in\gotA_{\sgotP}$  such that
$x\,y=x$  and  \[
\|x\,a\,\alpha(x)\|<\varepsilon\;\;\;.
\]
If $\gotA$ is simple and unital then $\alpha$ is properly outer if
and only if $\alpha$ is outer.

\begin{thm}\label{tinvst2.1}
Let $G$ be an amenable discrete group, $\alpha$ an action of $G$ on a
$C^*$-algebra $\gotA$ such that $\alpha_g$ is properly outer for
$g\in G\backslash\{e\}$ and assume there exists a faithful $\alpha$-covariant
irreducible representation of $\gotA$.
Then for any $\alpha$-invariant state $\omega$ on $\gotA$, any finite subset
$F\subseteq G$, any finite subset $\gotF\subseteq\gotA$ and any $\varepsilon>0$
there exists a pure state $\varphi$ on $\gotA$ such that
\[
\|\varphi\circ\alpha_g-\varphi\|<\varepsilon
\]
for all $g\in F$ and
\[
|\varphi(x)-\omega(x)|<\varepsilon
\]
for all $x\in\gotF$.
Moreover, $\varphi$ can be taken to be a vector state in the given
$\alpha$-covariant representation.
\end{thm}

\begin{remarkn}\label{rinvst2.1} If $G$ is abelian and countable
and $\gotA$ is separable and prime then the existence of a faithful
$\alpha$-covariant irreducible representation is automatic (see
Theorem~3.3 in \cite{Kis4}).
\end{remarkn}

\begin{remarkn}\label{rinvst2.2} The hypotheses of the theorem do
not imply that the set of  $\alpha$-invariant pure states is
weak$^*$-dense in the set of $\alpha$-invariant states.
The set may be weak$^*$-closed (see Example \ref{einvst2.1}) and
there are even examples satisfying the hypotheses of the theorem
for which there are no $\alpha$-invariant pure states (see Examples
\ref{einvst2.2} and \ref{einvst2.3}).
\end{remarkn}

As a preliminary to the proof of the theorem we establish two lemmas.
\begin{lemma}\label{linvst2.1}
Let $\gotA$ be a $C^*$-algebra, $\varphi$ a pure state on $\gotA$
such that the associated representation is faithful, $\varepsilon>0$
and $\gotF\subset\gotA$ a finite subset of the algebra.
Then there exists an element $e\in\gotA_{\sgotP}$ such that
$\varphi(e)=1$ and
\[
\|e\,x\,e-\varphi(x)\,e^2\|<\varepsilon
\]
for all $x\in\gotF$.
\end{lemma}
\proof\
By Proposition 3.13.6 in \cite{Ped} the support projection $p$ of
$\varphi$ in the bidual $\gotA^{**}$ of $\gotA$ is closed.
If, temporarily, we assume that $\gotA$ is separable then there is
a decreasing sequence $z_n$ of positive elements in $\gotA$ such
that $z_n\searrow p$.
Setting $z=\sum_n2^{-n}z_n$ and taking suitable
functions of $z$ we can construct a
sequence $\{e_n\}$ of positive elements such that $0\leq e_n\leq1$,
\[
e_n\,e_m=e_{\max\{n,m\}}
\]
for all $n,m\in\Ni$ and $e_n\searrow p$.
Thus if $p_n\in\gotA^{**}$ is the eigenprojection of $e_n$ corresponding to
eigenvalue one and $m>n$ then $e_m\,p_n=p_n\,e_m=e_m$.
We also have $p_n\searrow p$.
We may assume that all elements in $\gotF$ are selfadjoint.
For $x\in\gotF$ let $\sigma(p_n\,x\, p_n)$ denote the spectrum of
$p_n\,x\,p_n$ as an element in $p_n\,\gotA^{**}\,p_n$.
Then $\max \sigma(p_n\,x\,p_n)$ is a decreasing sequence, with limit $a$, and
$\min\sigma(p_n\,x\,p_n)$ is an increasing sequence, with limit $b\leq a$.
We argue that $b=a=\varphi(x)$.
If not there are states $\omega_a$ and $\omega_b$ on $\gotA$ with
$\omega_a(x)=a$,
 $\omega_a(e_k)=1$, $\omega_b(x)=b$,
 $\omega_b(e_k)=1$ for all $k$.
But then $\supp \omega_a\leq\lim_n e_n=p$ and $\supp \omega_b\leq p$.
Since $p$ is a one-dimensional projection in $\gotA^{**}$ by purity of $\varphi$
one then has $\omega_a=\omega_b=\varphi$ and $a=b=\varphi(x)$.
Also, as $\max \sigma(p_n\,x\,p_n)\to\varphi(x)$ and
$\min \sigma(p_n\,x\,p_n)\to\varphi(x)$ one has
$\|p_n\,x\,p_n-\varphi(x)\,p_n\|\to 0$.
Thus given $\varepsilon>0$ and $\gotF$ we can find an $n$ with
\[
\|p_n\,x\,p_n-\varphi(x)\,p_n\|<\varepsilon
\]
for all $x\in\gotF$.
But if $m>n$ then
\begin{eqnarray*}
\|e_m\,x\,e_m-\varphi(x)\,e_m^2\|&=&
\|e_m\,p_n\,x\,p_n\,e_m-\varphi(x)\,e_m\,p_n\,e_m\|\\[5pt]
&\leq&\|p_n\,x\,p_n-\varphi(x)\,p_n\|<\varepsilon
\end{eqnarray*}
and the lemma is established for separable $\gotA$.

If $\gotA$ is nonseparable the proof of the lemma can be reduced to the
separable
case by the following argument.

Adopt the hypotheses of the lemma and let $\gotB_1$ be the separable
$C^*$-subalgebra of $\gotA$ generated by $\gotF$.
We next construct inductively a sequence $\gotB_1\subseteq
\gotB_2\subseteq\ldots$
of separable $C^*$-subalgebras of $\gotA$ as follows.
Assume $\gotB_n$ has been constructed and let $\{x_k\}$ be a dense
sequence in $\gotB_n$.
If $\Phi$ is the cyclic vector corresponding to $\varphi$ then $\{x_k\Phi\}$ is
dense in $[\,\gotB_k\Phi\,]$.
But if $k,m$ are such that
\[
|\,\|x_k\Phi\|-\|x_m\Phi\|\,|<n^{-1}\|x_m\Phi\|
\]
then by Kadison's transitivity theorem there is a
$u_{k;m}\in\gotA$ such that $\|u_{k;m}\|\leq1$ and
\[
\|x_k\Phi-u_{k;m}x_m\Phi\|<n^{-1}\|x_m\Phi\|\;\;\;.
\]
Also, for each $m$, there is a  $u_m\in\gotA$ such that
$\|u_m\|\leq1$ and
\[
\varphi(u_m\,x_m\,x_m^*\,u_m^*)\geq(1-n^{-1})\|x_m\,x_m^*\|\;\;\;.
\]
Now let $\gotB_{n+1}$ be the $C^*$-algebra generated by  $\{x_m\}, \{u_{k;m}\},
\{u_m\}$.

Set $\gotB=\overline{\bigcup_n\gotB_n}$.
Then $\gotB$ is a separable $C^*$-subalgebra of $\gotA$ containing $\gotF$
and having
the property that if $\psi_1, \psi_2\in[\,\gotB\Phi\,]$  with
$\|\psi_1\|=\|\psi_2\|\neq0$ and $\varepsilon>0$ then there is a
$u\in\gotB$ such that
\[
\|\psi_2-u\psi_1\|\leq\varepsilon\,\|\psi_1\|
\]
and it follows that the representation of $\gotB$ on $[\,\gotB\Phi\,]$ is
irreducible.
Furthermore if $x\in\gotB$ and $x\neq0$ there is a
$u\in\gotB$ with $\varphi(u\,x\,x^*\,u^*)>0$ and hence this representation
is faithful. Thus we may apply the separable case of the lemma, which has
already been established, to $\gotB$, ${\varphi}|_{\sgotB}$ and $\gotF$
instead of $\gotA$, $\varphi$ and $\gotF$ to conclude the general validity
of the lemma. \hfill\Box

\begin{lemma}\label{linvst2.2}
Let $G$ be a discrete group  and $\alpha$ an action of $G$ on a prime
$C^*$-algebra such that $\alpha_g$ is properly outer for $g\in
G\backslash\{e\}$.
Let $\omega$ be an $\alpha$-invariant state.
For all finite subsets $\Lambda\subseteq G$, all $\delta>0$ and all finite
subsets $\gotF\subseteq\gotA$ there exists an element
$p\in\gotA_{\sgotP}$ with
\begin{description}
\item[{\rm(i)}] $\;\;\;$ if $g\neq h$, $g,h\in\Lambda$ and
$x\in\gotF\bigcup\{I\}$ then $\|\alpha_g(p)\,x\,\alpha_h(p)\|<\delta$
\item[{\rm(ii)}] $\;\;$ if  $g\in\Lambda$ and $x\in\gotF$ then
$\|p\,\alpha_{g^{-1}}(x)\,p-\omega(x)\,p^2\|<\delta$.
\end{description}
\end{lemma}
\proof\
The proof closely follows an argument of \cite{Kis5}.
First, by replacing $G$ by the group generated by $\Lambda$ and using an
argument similar to that occurring at the end of the proof of
Lemma \ref{linvst2.1} one may assume that the algebra $\gotA$
is separable.
The argument in the present setting is, however, more complicated
since one now has to ensure that the restriction of $\alpha$ to
the new subalgebra is still properly outer.
By Glimm's theorem (see, for example, \cite{Gli}), there is for each
$\delta>0$ a pure state $\varphi$ on $\gotA$ defining a faithful
representation such that
\[
|\varphi(\alpha_g(x))-\omega(\alpha_g(x))|=
|\varphi(\alpha_g(x))-\omega(x)|<\delta
\]
for any $x\in\gotF$ and  $g\in\Lambda^{-1}$.
But by Lemma \ref{linvst2.1} there is an element $e\in\gotA_{\sgotP}$
such that $\varphi(e)=1$ and
\[
\|e\,\alpha_g(x)\,e-\varphi(\alpha_g(x)\, e^2\|<\delta
\]
for all $x\in\gotF$ and $g\in\Lambda^{-1}$.
Therefore, slightly modifying $e$ by spectral theory, we may assume
that there is a  $q'\in\gotA_{\sgotP}$ with $\varphi(q')=1$ and
$\|q'\|=1$  such that  $q'\,e=q'$.
Then the above estimate is  still valid with $e$ replaced by any such
$q'$.

Since $\alpha_g$ is properly outer for each $g\in G $
there is a $p\in\gotA_{\sgotP}$ with $e\,p=p$  such that
\[
\|\alpha_g(p)\,x\,\alpha_h(p)\|=
\|\alpha_{h^{-1}g}(p)\,\alpha_{h^{-1}}(x)\,p\|<\delta
\]
for all $g,h\in\Lambda$ with $g\neq h$ and all $x\in\gotF\bigcup \{I\}$.
For a single choice of $g,h$ and $x$ the existence of this $p$ follows
from proper outerness of $\alpha_{h^{-1}g}$ by Lemma 1.1 of \cite{Kis4}
but going through the finite list of $g,h$ and $x$ one successively
constructs new $p'$ with $p'\,p=p'$ satisfying the estimate above
with the new $p'$ for the new triple $(g,h,x)$.
Then the estimate holds with $p$ replaced by $p'$ for the earlier
elements in the list.
Continuing with this process one obtains the sought after $p$ and
property (i) is immediate.
As $p\,e=p$ one has
\[
\|p\,\alpha_{g^{-1}}(x)\,p-\varphi(\alpha_{g^{-1}}(x))\,p^2\|<
\delta
\]
for all $x\in\gotF$ and $g\in\Lambda$.
But as $\varphi(\alpha_{g^{-1}}(x))$ is close to $\omega(x)$ one
finds
\[
\|p\,\alpha_{g^{-1}}(x)\,p-\omega(x)\,p^2\|<2\,\delta
\]
and so (ii) is valid.\hfill\Box

\bigskip

\noindent{\bf Proof of Theorem \ref{tinvst2.1}}$\;$ Let $\pi$ be a
faithful $\alpha$-covariant irreducible representation of $\gotA$ on a
Hilbert space $\ch$ and choose a finite subset $\Lambda\subseteq G$
such that $F\subseteq \Lambda$ and
\[
|\,\Lambda\,\triangle\,h\Lambda\,|/|\,\Lambda\,|<\delta
\]
for $h\in F$.
This is possible by the amenability of $G$.
Next choose $p$ as in Lemma \ref{linvst2.2} and $q\in\gotA_{\sgotP}$
with $p\,q=q$ where the $\Lambda$ in the lemma is replaced by
$F\Lambda\bigcup\,\Lambda$ and the $\delta$ by
$\delta/|\Lambda|^2$.
Then choose a unit vector $\psi\in\pi(q)\ch$,
i.e., $\psi=\pi(q)\xi$ for some $\xi\in\ch$.
Therefore
\[
\pi(p)\psi=\pi(p\,q)\xi=\pi(q)\xi=\psi\;\;\;.
\]
Now define
\[
\Phi=|\Lambda|^{-1/2}\sum_{g\in\Lambda}U_g\psi
\]
where $U$ is the unitary representation of $G$ on $\ch$ such that
\[
U_g\pi(x)U^*_g=\pi(\alpha_g(x))
\]
for all $x\in\gotA$ and $g\in G$.
Next define
\[
\varphi(x)=(\pi(x)\Phi,\Phi)/\|\Phi\|^2
\]
for all $x\in\gotA$.
Then $\varphi$ is a pure state on $\gotA$ and we will verify that
$\varphi$ has the properties claimed in Theorem \ref{tinvst2.1}.

First, if $h\in F$ one has
\[
\|U_h\Phi-\Phi\|=
|\Lambda|^{-1/2}\|\sum_{g\in\Lambda}(U_{hg}\psi-U_g\psi)\|
\;\;\;.
\]
In the latter sum all $g\in\Lambda\bigcap h^{-1}\Lambda$ cancel and
the remainder is a surface term
\[
|\Lambda|^{-1/2}
\|\sum_{g\in h\Lambda\backslash\Lambda}
U_g\psi-\sum_{g\in \Lambda\backslash h\Lambda}
U_g\psi\|
\;\;\;.
\]
But if $g_1, g_2\in F\Lambda\bigcup\Lambda$ with $g_1\neq g_2$ then
\begin{eqnarray*}
|(U_{g_1}\psi,U_{g_2}\psi)|&=&
|(U_{g_1}\pi(p)\psi,U_{g_2}\pi(p)\psi)|\\[5pt]
&=&
|(\pi(\alpha_{g_1}(p))U_{g_1}\psi,\pi(\alpha_{g_2}(p))U_{g_2}\psi)|\\[5pt]
&\leq&
\|\alpha_{g_1}(p)\alpha_{g_2}(p)\|<\delta/|\Lambda|^2
\end{eqnarray*}
by Lemma \ref{linvst2.2}.
Thus
\begin{eqnarray*}
\|U_h\Phi-\Phi\|^2&\leq&|\Lambda|^{-1}
\sum_{g_1,g_2\in h\Lambda\triangle\Lambda}
|(U_{g_1}\psi,U_{g_2}\psi)|\\[5pt]
&=&|\Lambda|^{-1}\sum_{g\in h\Lambda\triangle\Lambda}1
+|\Lambda|^{-1}
\sum_{{g_1\neq g_2}\atop{g_1,g_2\in h\Lambda\triangle\Lambda}}
|(U_{g_1}\psi,U_{g_2}\psi)|\\[5pt]
&\leq&|h\Lambda\triangle\Lambda|/|\Lambda|+
\delta|h\Lambda\triangle\Lambda|^2/|\Lambda|^2\\[5pt]
&\leq&\delta+4\delta=5\delta\;\;\;.
\end{eqnarray*}
In addition
\begin{eqnarray*}
|\,\|\Phi\|^2-1\,|&\leq&
|\Lambda|^{-1}\sum_{{g_1,g_2\in\Lambda}\atop{g_1\neq g_2}}
|(U_{g_1}\psi,U_{g_2}\psi)|\\
&\leq&\delta\,|\Lambda|^{-3}|\Lambda|^2\leq\delta\,|\Lambda|^{-1}
\;\;\;.
\end{eqnarray*}
It follows that
\[
\|\varphi\circ\alpha_h-\varphi\|\leq2(5\delta)^{1/2}(1+\delta)
\]
for all $h\in F$ so the first estimate of Theorem \ref{tinvst2.1} is
valid.
Secondly, one checks the second estimate as follows.
One has
\begin{eqnarray*}
\|\Phi\|^2\varphi(x)&=&(\pi(x)\Phi,\Phi)\\[5pt]
&=&|\Lambda|^{-1}\sum_{g,h\in\Lambda}
(\pi(x\,\alpha_g(p))U_g\psi,\pi(\alpha_h(p))U_h\psi)\\[5pt]
&=&|\Lambda|^{-1}\sum_{g,h\in\Lambda}
(\pi(\alpha_h(p)\,x\,\alpha_g(p))U_g\psi,U_h\psi)
\;\;\;.
\end{eqnarray*}
But one also has
\begin{eqnarray*}
\|\Phi\|^2\omega(x)&=&(\omega(x)\Phi,\Phi)\\[5pt]
&=&|\Lambda|^{-1}\sum_{g,h\in\Lambda}
(\pi(\alpha_h(p)\,\omega(x)\,\alpha_g(p))U_g\psi,U_h\psi)
\end{eqnarray*}
and therefore
\[
\|\Phi\|^2|\varphi(x)-\omega(x)|\leq
|\Lambda|^{-1}\sum_{g,h\in\Lambda}
\|\alpha_h(p)\,x\,\alpha_g(p)-\omega(x)\,\alpha_h(p)\,\alpha_g(p)\|
\;\;\;.
\]
If, however, $x\in\gotF$ and $g\neq h$ it follows from Lemma
\ref{linvst2.2} (ii) that
\begin{eqnarray*}
\|\alpha_h(p)\,x\,\alpha_g(p)-\omega(x)\,\alpha_h(p)\,\alpha_g(p)\|
&\leq&\|\alpha_h(p)\,x\,\alpha_g(p)\|+
\|x\|\cdot\|\alpha_h(p)\,\alpha_g(p)\|\\[5pt]
&\leq&\delta(1+|\gotF|)/|\Lambda|^2
\end{eqnarray*}
where $|\gotF|=\max\{\|x\|:x\in\gotF\}$ and if $g=h$ it follows from
Lemma \ref{linvst2.2} (iii) that
\[
\|\alpha_g(p)\,x\,\alpha_g(p)-\omega(x)\,\alpha_g(p)^2\|
=\|p\,\alpha_{g^{-1}}(x)\,p-\omega(x)\,p^2\|\leq\delta/|\Lambda|^2
\;\;\;.
\]
Hence
\[
|\varphi(x)-\omega(x)|\leq|\Lambda|^2\,\Big(\delta(1+|\gotF|)/|\Lambda|^2
+|\Lambda|\,\delta/|\Lambda|^2\Big)(1+\delta)
\]
and the second estimate of Theorem \ref{tinvst2.1} follows.\hfill\Box

\medskip

\begin{exam}\label{einvst2.1} It is not true in general that for
$\gotA, G$ and $\alpha$ satisfying the hypotheses of Theorem
\ref{tinvst2.1}, even with $G=\Zi$, that the set of $G$-invariant
pure states on $\gotA$ is weakly$^*$-dense in the set of $G$-invariant
states.
Let $\gotA$ be the UHF-algebra $\gotA=\otimes^\infty_{n=1}M_2$ and let
\[
\alpha=\bigotimes^\infty_{n=1}{\rm Ad}
\left(
\begin{array}{cc}
1  &  0\\
0  &e^{2\pi i\theta_n}
\end{array}
\right)
\]
where the irrational numbers $\theta_i$ are rationally independent.
This action is almost periodic with fixed point algebra equal to the
maximal abelian subalgebra $\otimes^\infty_{n=1}\Ci^2$.
Thus invariant states correspond to the probability measures on
$\times^\infty_{n=1}\{0,1\}$ and extremal invariant states to point
evaluations on the Cantor set. So the extremal invariant states are
the pure product measures
\[
\rho(x)=\Tr\left(x\bigotimes^\infty_{n=1}
\left(
\begin{array}{cc}
\rho_n  &  0\\
0  &1-\rho_n
\end{array}
\right)\right)
\]
on $\gotA$ where $\rho_n\in\{0,1\}$.
Hence all extremal invariant states are pure and the set of extremal
invariant states is closed in the weak$^*$-topology.
\end{exam}

\medskip

\begin{exam}\label{einvst2.2}
We next construct a $C^*$-dynamical system $(\gotA, \alpha,\Zi)$
satisfying the hypotheses of Theorem \ref{tinvst2.1} for which
there are no invariant pure states!
Let $\theta\in[0,1]$ be an irrational number and let $\gotA_\theta$
be the universal $C^*$-algebra generated by two unitaries $U$, $V$
with $VU=e^{2\pi i\theta}UV$, \cite{Rie}. The algebra
$\gotA_\theta$ is simple and carries an action $\beta$ of the
two-torus given by
\[
\beta_{(z_1,z_2)}(U)=z_1U\;\;\;\;\;,\;\;\;\;\;
\beta_{(z_1,z_2)}(V)=z_2V
\]
for $z_1,z_2\in\Ti$.
This action is ergodic and the only invariant state for the action
$\beta$ is the trace state $\omega$ defined by
\[
\omega(U^nV^m)=\left\{\begin{array}{ll}
                          0 &\mbox{if $(n,m)\neq(0,0)$}\\
                          1 &\mbox{ if $n=m=0$.}
                                          \end{array}\right.
\]
Note that
\begin{eqnarray*}
U^nV^mU(U^nV^m)^*&=&e^{2\pi im\theta}U\\[5pt]
U^nV^mV(U^nV^m)^*&=&e^{-2\pi in\theta}V
\end{eqnarray*}
and it can be proved that $\beta_{(z_1,z_2)}$ is inner if and only if
\[
z_1,z_2\in(\theta\Zi+\Zi)/\Zi\subseteq \Ri/\Zi=\Ti
\]
(see \cite{ElR}, \cite{Kis7}, \cite{BKR2}).
(One well-known argument for this is the following: If $\beta_{(z_1,z_2)}$
is implemented by a unitary $W\in\gotA$ then $WUW^*=z_1U$ and
$WVW^*=z_2V$.
Hence $UWU^*={\overline{z_1}}W$ and $VWV^*={\overline{z_2}}V$.
Thus $W$ is an eigenunitary for all the automorphisms $\Ad(U^nV^m)$.
But as these are dense in $\beta_{\Ti^2}$ it follows that $W$ is an
eigenunitary
for the action $\beta$.
Hence $W=\lambda U^nV^m$ for some $\lambda\in\Ti$ and $n,m\in\Zi$.
Or just note that  as $\beta_{(z_1,z_2)}$ commutes
with all $\beta_{(w_1,w_2)}$,  one has $\Ad \beta_w(W) = \Ad W$ for all $w \in
\Ti^2$, and hence W is an eigenelement for the action $\beta$ by simplicity of
$\gotA$.)
Now choose irrational numbers $\theta_1,\theta_2$ such that
$1,\theta, \theta_1$ and $\theta_2$ are rationally independent, put
\[
g=(e^{2\pi i\theta_1},e^{2\pi i\theta_2})\in\Ti^2
\]
and $\alpha=\beta_g$.
Since $\{\,g\,|\,n\in\Zi\,\}$ is dense in $\Ti^2$ the only
$\alpha$-invariant state is the trace state $\omega$.
Also, since
\[
\{\,g^n\,|\,n\in\Zi\,\}\bigcap \,
(e^{2\pi i\theta\Zi},e^{2\pi i\theta\Zi})=0\;\;\;,
\]
all the automorphisms $\alpha^n$, $n\in\Zi\backslash\{0\}$ are outer.
(A simple direct argument for this is the following:
If $\alpha^n=\Ad(U)$ for some $ n\in\Zi\backslash\{0\}$ and $U\in\gotA$
then any character on $C^*(U)$ extends to a pure $\alpha^n$-invariant state
on $\gotA$.
But the powers $(\alpha^n)^m$, $m\in\Zi$, are dense in $\beta_{\Ti^2}$ so
this state must be the trace.
This contradiction shows that all  $\alpha^n$, $n\in\Zi\backslash\{0\}$
are outer.)
Furthermore, $(\gotA, \alpha, \Zi)$ admits a faithful
$\alpha$-covariant irreducible representation by
Remark~\ref{rinvst2.1}.
Hence all the hypotheses of Theorem~\ref{tinvst2.1} are satisfied
but there are no $\alpha$-invariant pure states on $\gotA$.
\end{exam}

\medskip

\begin{exam}\label{einvst2.3}
There exists a UHF-algebra $\gotA$ and an automorphism $\alpha$
of $\gotA$ such that the associated action of $\Zi$ satisfies all
the assumptions of Theorem~\ref{tinvst2.1} but nevertheless there
is no invariant pure state.
Take $p,q\in \{2,3,\ldots\}$ and set
$\gotA=M_{p^\infty}\otimes M_{q^\infty}$,
$\alpha=\alpha_1\otimes\alpha_2$ where $\alpha_1$ is the shift on
$M_{p^\infty}=\otimes^\infty_{-\infty}M_p$ and $\alpha_2$ is the
automorphism on $M_{q^\infty}$ constructed by Connes in \cite{Con}
with outer invariant $(q,\gamma)$ where $\gamma\neq1$ is a $q$-th
root of unity.
This means that $q$ is the smallest natural number such that
$\alpha_2^q$ is inner in the trace representation of $M_{q^\infty}$
and if $\alpha_2^q=\Ad(u)$ for a unitary $u$ in the weak closure then
$\alpha_2(u) =\gamma u$. (We have also used $\alpha_2$ to  denote
the extension of the automorphism to the weak closure.)
Explicitly, $\alpha_2$ is defined as follows.
Set $M_{q^\infty}=\otimes^\infty_{m=1}M_q^{(m)}$ where
$M_q^{(m)}\cong M_q$.
Then if $\theta$ is the one-sided shift on $M_{q^\infty}$, define
the unitaries $U\in M^{(1)}_q$ and $v\in M^{(1)}_q\otimes M^{(2)}_q$
by
\[
U= \left ( \begin{array}{cccccc}
          \gamma & 0 &  \cdots & 0 & 0  \\
        0 & \gamma^2 & \cdots & 0 & 0  \\
\vdots & \vdots  & \ddots & \vdots & \vdots \\
      0 & 0 & \cdots & \gamma^{q-1} & 0  \\
         0 & 0  & \cdots & 0 & \gamma^q
\end{array} \right ) \;
\;\;\;{\rm and} \;\;\;\;
v= \left ( \begin{array}{ccccc}
 \,0\, & \,1\, &  \,\cdots\, & \,0\, & \,0\,  \\
          0 & 0 &  \ddots & 0 & 0  \\
\vdots & \vdots  & {} & \ddots & \vdots \\
          0 & 0 &  \cdots & 0 & 1  \\
          \theta(U^*) & 0 &  \cdots & 0 & 0
\end{array} \right )
\]
and, following \cite{Con}, set
\[
\alpha_2=\lim_{n\to\infty}
\Ad(v\,\theta(v)\,\theta^2(v)\ldots\theta^n(v))\;\;\;.
\]
Then it follows that $\alpha_2^q= \Ad U$ and $\alpha_2(U)=\gamma U$
by Proposition~1.6 in \cite{Con}.
Now if $\varphi$ is an $\alpha$-invariant pure state on $\gotA$ it
follows from the fact that $\alpha_1$ is asymptotically abelian that
$\varphi=\varphi_1\otimes\varphi_2$ where $\varphi_i$ is an $\alpha_i$-invariant
pure state for $i=1,2$.
But as $\gamma\neq0$ it follows that $\alpha_2$ is not a
quasi-product action of the group $\Zi_{rq}$, where $r\in\Ni$ is the
smallest number with $\gamma^r=1$.
This can be seen from the remark prior to Theorem~\ref{tinvst3.1}
below and from this remark it also follows that $(M_{q^\infty},
\alpha_2)$ does not admit an invariant pure state.
More specifically, if there exists an $\alpha_2$-covariant irreducible
representation, there is a unitary operator $V$ on the representation space with
$\alpha_2=\Ad V$, thus $\alpha_2^q= \Ad V^q=Ad U$ and $V^q=\eta U$ for an
$\eta\in\Ti$ by irreducibility.
But then $\alpha_2(V)=VVV^*=V$, so $\alpha_2(V^q)=V^q$ and
$\alpha_2(U)=U\neq\gamma U$, where $\alpha_2$ also denotes the extension of
$\alpha_2$ to the weak closure.
This contradiction establishes that $(M_{q^\infty},
\alpha_2)$ does not admit a covariant irreducible representation and thus
it does
not admit an invariant pure state.
Furthermore, as
$\alpha_1$ is asymptotically abelian, $\alpha_1^n$ is outer for all
$n\in\Zi\backslash\{0\}$. Thus $(\gotA, \alpha, \Zi)$ satisfies all
the conditions of Theorem~\ref{tinvst2.1} but does not admit an
$\alpha$-invariant pure state.
\end{exam}

\section{Density of invariant pure states for compact group actions}

Let $\alpha$ be an action of a compact group $G$ on a $C^*$-algebra
$\gotA$ and let
\[
x\in\gotA\mapsto E(x)=\int_Gdg\,\alpha_g(x)\in\gotA^\alpha
\]
be the associated projection onto the fixed point algebra
$\gotA^\alpha$.
If $\omega$ is a state on $\gotA^\alpha$ then $\omega\circ E$ is an
$\alpha$-invariant state on $\gotA$ and the map
$\omega\mapsto\omega\circ E$ is an affine homeomorphism between the
compact convex set $E_{\gotA^\alpha}$ of states on $\gotA^\alpha$ and
the compact convex set $E^G_\alpha$ of $\alpha$-invariant states on
$\gotA$.
Thus pure states on $\gotA^\alpha$ correspond exactly to extremal
invariant states on $\gotA$.
In this section we address the question of when purity of $\omega$
on $\gotA^\alpha$ implies purity of $\omega\circ E$ on $\gotA$.
Simple examples show that it may happen that no extremal invariant
states are pure, e.g., if $\gotA$ is abelian and there are no fixed
points in the spectrum of $\gotA$ for the $\alpha$ action.
More interesting examples are the canonical actions of the
two-torus $\Ti^2$  on the irrational rotation algebras mentioned in
Example~\ref{einvst2.2}.
These actions are ergodic and the only
invariant state is the trace state.
But in this example the automorphisms are inner for a dense set of
points in $\Ti^2$.
The next theorem demonstrates that this situation does not
generally arise in the converse case of quasi-product actions.

Recall from \cite{BEK}, \cite{BEEK} that a faithful action
$\alpha$ of a compact group $G\neq0$ on a $C^*$-algebra $\gotA$  is
called a quasi-product action if it satisfies any of the following
equivalent conditions;

\smallskip

\noindent (1)$\;$ If $x,y\in \gotA$ then
$\sup\{\,\|x\,a\,y\|\,|\,a\in\gotA^\alpha\,,\,\|a\|=1\,\}
=\|x\|\cdot\|y\|$.

\smallskip

\noindent (2)$\;$ There exists a faithful irreducible
representation $\pi$ of $\gotA$  such that $\pi|_{\sgotA^\alpha}$ is
irreducible.

\smallskip

\noindent (3)$\;$ There exists an $\alpha$-invariant pure state
$\omega$ on $\gotA$ such that $\pi_{\omega|_{\sgotA^\alpha}}$  is a
faithful representation for $\gotA^\alpha$.

\smallskip

\noindent (4)$\;$ For any sequence $(\zeta_n)$ of
finite-dimensional unitary representations of $G$ there exists an
$\alpha$-invariant sub-$C^*$-algebra $\gotB$ of $\gotA$  and a
closed $\alpha^{**}$-invariant projection $q$ in the bidual
$\gotA^{**}$ of $\gotA$ such that
\[
q\in\gotB'\;\;\;\;\;,\;\;\;\;\;q\,\gotA\, q=\gotB q\;\;\;,
\]
 $q\in\gotI^{**}\subseteq\gotA^{**}$ for any non-zero two-sided
ideal $\gotI$ of $\gotA$ and the $C^*$-dynamical system
$(\gotB q, G, \alpha^{**}|_{\sgotB q})$ is isomorphic to the
product system
\[
(\otimes^\infty_{n=1}M_{\dim(\zeta_n)}, G,
\otimes^\infty_{n=1}\Ad\zeta_n)\;\;\;.
\]

If $G$ is abelian these conditions are also equivalent to

\smallskip

\noindent (5)$\;$ $\gotA$ and $\gotA^\alpha$ are prime and
$\alpha_g$ is properly outer for each $g\in G\backslash\{0\}$.

\smallskip

It is maybe now not too surprising that these conditions are also
equivalent with the density of the invariant pure states in the set
of invariant states.
\begin{thm}\label{tinvst3.1}
Let $\alpha$ be a faithful action of a compact group $G\neq0$ on a
separable  $C^*$-algebra $\gotA$.
The following conditions are equivalent:
\begin{description}
\item[{\rm(1)}]
$\alpha$ is a quasi-product action,
\item[{\rm(2)}]
The set of  $\alpha$-invariant pure states on $\gotA$ is
weak$^*$-dense in the set of $\alpha$-invariant states.
\end{description}
\end{thm}
In proving Theorem~\ref{tinvst3.1} we begin by noting that
(2)$\Rightarrow$(1) follows straightforwardly from known results
as follows.

First, the characterization (3) of quasi-product states is
equivalent to

\medskip

\noindent($3'$)$\;$ $\gotA^\alpha$ is prime and there exists a family
$\{\omega_i\}$ of $\alpha$-invariant pure states $\omega_i$ on
$\gotA$ such that $\oplus_i\pi_{\omega_i|_{\sgotA^\alpha}}$ is a
faithful representation of $\gotA^\alpha$.

\medskip

Clearly (3)$\Rightarrow$($3'$) and the proof of the converse is
indicated in \cite{Kis6}.
Briefly one has to check the proof of Theorem~3.1 in \cite{BKR}.
But referring to Condition~(6) in \cite{BEK}, which is the same as
Condition~(4) above, one has to prove in addition that one may
choose $q\in\gotI^{**}$ for any non-zero two-sided ideal $\gotI$ of
$\gotA$.
This, however, follows once one can show from ($3'$) that $\gotA$ is
prime.
So suppose {\it ad absurdum} that $\gotA$ has two ideals $\gotI_1$,
$\gotI_2$ with zero intersection.
Then there is an invariant pure state $\omega$ such that
$\omega|_{\sgotI_1}$ is a state and $\omega(\gotI_2)=0$.
Since $\omega$ is also zero on the invariant ideal $\gotK_2$
generated by $\gotI_2$ one then has $\gotK_2\gotI_1=0$.
This remains true if $\gotI_1$ is replaced by the invariant ideal
generated by $\gotI_1$.
But intersecting with $\gotA^\alpha$ one obtains two non-zero
ideals in $\gotA^\alpha$ with zero intersection which contradicts
the assumption that $\gotA^\alpha$ is prime.
Thus $\gotA$ is prime.

Returning to the proof of (2)$\Rightarrow$(1) in
Theorem~\ref{tinvst3.1} it now suffices to show that (2) implies
that $\gotA^\alpha$ is prime, since the other conditions in ($3'$)
obviously follow from (2).
But (2) clearly implies that the extremal $\alpha$-invariant states
are dense in the set of $\alpha$-invariant states and, by the
initial remarks in this section, this means that the pure states on
$\gotA^\alpha$ are dense in the states on $\gotA^\alpha$.
This means that $\gotA^\alpha$ is prime (and antiliminal).

In the proof of  the implication (1)$\Rightarrow$(2) in
Theorem~\ref{tinvst3.1} the characterization  (2) of quasi-product
actions is the most convenient.

We first focus on the case
that $G$ is abelian.
For $\gamma\in{\hat G}$ let
\[
\gotA^\alpha(\gamma)=\{x\in\gotA\,|\,\alpha_g(x)=\gamma(g)x\,\}
\]
be the corresponding spectral subspace.
If $\varphi$ is a state on $\gotA^\alpha$ and
$x\in\gotA^\alpha(\gamma)$ then $x^*x\in\gotA^\alpha$ and if
$\varphi(x^*x)=1$ then $y\mapsto \varphi(x^*yx)=(x\varphi x^*)(y)$
is a state on $\gotA^\alpha$.
In this situation let $\ch=\ch_{\varphi\circ E}$,
$\Omega=\Omega_{\varphi\circ E}$, $\pi=\pi_{\varphi\circ E}$ denote
the associated Hilbert space, cyclic vector and representation,
and $U$ the associated representation of $G$.
Thus $U_g=\sum_{\gamma\in{\hat G}}\gamma(g)P_\gamma$ where
$P_\gamma=[\,\pi(\gotA^\alpha(\gamma))\Omega\,]$.
If $P_\gamma\neq0$ let $\rho_\gamma$ be the representation of
$\gotA^\alpha$ obtained by restricting $\pi|_{\gotA^\alpha}$ to
$P_\gamma\ch=[\,\pi(\gotA^\alpha(\gamma))\Omega\,]$.

We can now formulate  general criteria for purity of an extremal
invariant state. Most of this proposition is more or less known.
\begin{prop}\label{pinvst3.1}
Let $\gotA$ be a $C^*$-algebra, $\alpha$ an action of a compact
abelian group on $\gotA$ and $\varphi$ a state on $\gotA^\alpha$.
The following conditions are equivalent:
\begin{description}
\item[{\rm(1)}]
$\varphi\circ E$ is a pure state on $\gotA$,
\item[{\rm(2)}]
$\varphi$ is a pure state on $\gotA^\alpha$ and
$P_0=[\,\pi(\gotA^\alpha)\Omega\,]\in\pi(\gotA^\alpha)''$,
\item[{\rm(3)}]
$\varphi$ is a pure state on $\gotA^\alpha$ and
$P_\Omega=[\,\Ci\Omega\,]\in\pi(\gotA^\alpha)''$,
\item[{\rm(4)}]
$\varphi$ is a pure state on $\gotA^\alpha$ and
$\rho_0\disjoint \rho_\gamma \;\;\;{\rm for}\;\gamma\neq0$,
\item[{\rm(5)}]
$\varphi$ is a pure state on $\gotA^\alpha$ and
for all $\gamma\in{\hat G}\backslash\{0\}$ such that
$\varphi(\gotA^\alpha(\gamma)^*\gotA^\alpha(\gamma))\neq0$ there is an
$x\in\gotA^\alpha(\gamma)$ such that $\varphi(x^*x)=1$ and
$x\varphi x^*\disjoint\varphi$.
\end{description}

Furthermore, all these conditions imply that the representations
$\rho_\gamma$ of $\gotA^\alpha$ are irreducible and mutually disjoint
whenever they are non-zero.
\end{prop}
\proof\
First assume Condition (1) is valid.
By covariance of $\pi$ one has
\[
\gotB(P_\gamma\ch)=P_\gamma\gotB(\ch)P_\gamma
=P_\gamma\pi(\gotA)''P_\gamma
=\pi(\gotA^\alpha)''P_\gamma
\]
so the representations $\rho_\gamma$ of $\gotA^\alpha$ are
irreducible, when they are non-zero.
But as $P_\gamma$ is $\Ad(U_g)$-invariant it follows from covariance
that $P_\gamma\in\pi(\gotA^\alpha)''$.
Then it follows that
\[
\pi(\gotA^\alpha)'\bigcap\pi(\gotA^\alpha)''=\{P_\gamma\,|\,\gamma\in{\hat
G}\,\}''\;\;\;.
\]
Hence all the representations $\rho_\gamma$ are irreducible and
mutually disjoint, when they are non-zero.
Therefore Condition (1) implies each of the Conditions (2)--(5) as
well as the concluding statement in the proposition.

The proof of the converse implications relies on the following
result.
\begin{lemma}\label{linvst3.1}
Let $\gotA$ be a $C^*$-algebra, $\alpha$ an action of a compact
abelian group $G$ on $\gotA$ and $\varphi$ a pure state on
$\gotA^\alpha$.
It follows that all the non-zero representations $\rho_\gamma$ are
irreducible.
\end{lemma}
\proof\
If $\gamma=\hat e$ then the corresponding representation $\rho_0$ is
irreducible by assumption.
Next note that the set of vectors of the form
$\pi(x)\Omega$ with $x\in\gotA^\alpha(\gamma)$ are norm-dense in
$P_\gamma\ch$, so if one can show that the corresponding vector
states on $\gotA^\alpha$ are pure then $\rho_\gamma$ is necessarily
irreducible.
To this end pick an $x\in\gotA^\alpha(\gamma)$ and let
$\psi$ be a positive functional on $\gotA^\alpha$ such that $x\varphi
x^*\geq\psi$.
Thus there exists a $T\in\rho_\gamma(\gotA^\alpha)'$ with support
$[\,\rho_\gamma(\gotA^\alpha)\pi(x)\Omega\,]$ such that $T\geq0$ and
\[
\psi(a)=(\pi(a)\pi(x)\Omega,T\pi(x)\Omega)
\]
for all $a\in\gotA^\alpha$ (see, for example, \cite{BR1}
Theorem ~2.3.19).
But if $y\in\gotA^\alpha(\gamma)$ then
\[
y^*x\varphi x^*y\geq y^*\psi y
\]
as functionals on $\gotA^\alpha$.
But as $y^*x\in\gotA^\alpha$ it follows that $y^*x\varphi x^*y$ is a
vector functional in the $\rho_0$-representation and thus
$y^*x\varphi x^*y$ is a multiple of a pure state on $\gotA^\alpha$,
because $\varphi$ is assumed to be a pure state.
Thus there exists a $\lambda>0$ with
\[
y^*\psi y=\lambda\, y^*x\varphi x^*y
\]
and then
\[
y^*y\psi yy^*=\lambda\, yy^*x\varphi x^*yy^*
\;\;\;.
\]
Hence
\[
(T\pi(yy^*ayy^*)\pi(x)\Omega,\pi(x)\Omega)=
\lambda\,(\pi(yy^*ayy^*)\pi(x)\Omega,\pi(x)\Omega)
\]
for any $a\in\gotA^\alpha$.
In principle the value of $\lambda$ could
depend on $y$ but by letting $yy^*$ run through an approximate
identity for the ideal  $\gotA^\alpha(\gamma)\gotA^\alpha(\gamma)^*$
in $\gotA^\alpha$, or  replacing $x$ by $f(xx^*)x$ where $f$ has a
plateau at the value one and then choosing $y$ with $yy^*x=x$, one
deduces that $\lambda$ is in fact independent of $y$.
Therefore one concludes that $T=\lambda I$ on
$[\,\rho_\gamma(\gotA^\alpha)\pi(x)\Omega\,]$. It then follows that
$\psi=\lambda\, x\varphi x^*$ and $x\varphi x^*$ is a multiple of a
pure state.  Thus $\rho_\gamma$ is irreducible.\hfill$\Box$

\bigskip

\noindent{\bf Proof of Proposition \ref{pinvst3.1} continued}$\;$
We will prove that the conditions of the proposition are also
equivalent to

\smallskip

\noindent $(6)\;\;$ The representations $\rho_\gamma$ of
$\gotA^\alpha$ are irreducible, whenever  they are non-zero, and
mutually disjoint.

\smallskip

We have already argued that $(1)\Rightarrow (6)$ and Condition (6) is
equivalent to
\[
\pi_{\varphi\circ E}(\gotA^\alpha)''\bigcap
\pi_{\varphi\circ E}(\gotA^\alpha)'=
\pi_{\varphi\circ E}(\gotA^\alpha)'=\{P_\gamma\,|\,\gamma\in\hat
G\,\}'' \;\;\;.
\]
Thus if (6) is satisfied
\[
\pi_{\varphi\circ E}(\gotA)'\subseteq
\{P_\gamma\,|\,\gamma\in\hat G\,\}''
\;\;\;.
\]
But
$P_\gamma=[\,\pi_{\varphi\circ E}(\gotA^\alpha(\gamma))\Omega\,]$
and hence
\[
\pi_{\varphi\circ E}(\gotA^\alpha)P_{0}=
P_{\gamma}\pi_{\varphi\circ E}(\gotA^\alpha(\gamma))
\;\;\;.
\]
If $P_{\gamma}$ is non-zero then both sides of
the last equation are non-zero.
It now follows that $P_\gamma\not\in\pi_{\varphi\circ E}(\gotA)'$ if
$P_\gamma\neq0,1$ and $\gamma\neq0$ and one easily extends this
argument to show that \[
\pi_{\varphi\circ E}(\gotA)'=\pi_{\varphi\circ
E}(\gotA)'\bigcap\{P_\gamma|\,\gamma\in\hat G\,\}''=\Ci I
\;\;\;.
\]
(For example, the projections in $\{P_\gamma|\,\gamma\in\hat G\,\}$
are sums of the $P_\gamma$.) Thus $\varphi\circ E$ is pure and
$(6)\Rightarrow(1)$.

In order to show that each of the Conditions~(2)--(5)
implies Condition~(1) it therefore suffices to show that they imply
Condition~(6).
But by Lemma \ref{linvst3.1} each of these conditions imply that
$\rho_\gamma$ is zero or irreducible for each $\gamma\in\hat G$
and hence $(2)\Leftrightarrow(4)\Leftrightarrow(5)$.
Furthermore Condition~(3) implies that
\[
[\,\pi_{\varphi\circ E}(\gotA^\alpha)P_\Omega
\ch_{\varphi\circ E}\,]=
P_0\in\pi_{\varphi\circ E}(\gotA^\alpha)''
\]
so $(3)\Rightarrow(2)$.
Thus it remains to prove that $(4)\Rightarrow(6)$

Assume Condition~(4) is valid and that $\gamma_1$ and $\gamma_2$ are
such that $P_{\gamma_1}\neq0$ and $P_{\gamma_2}\neq0$.
Then the two representations $\rho_{\gamma_1}$ and $\rho_{\gamma_2}$
are irreducible and therefore either unitarily equivalent or
disjoint.
But if they are unitarily equivalent then for any
$x_i\in\gotA^\alpha(\gamma_i)$, $i=1,2$, with  $\varphi(x_i^*x_i)=1$,
the states $x_i\varphi x^*_i$ on $\gotA^\alpha$ are equivalent.
But by Kadison's transitivity theorem there then exists a unitary
$u\in\gotA^\alpha$ with
\[
x_2\varphi x_2^*=ux_1\varphi x_1^*u^*
\]
but then
\[
x_2^*x_2\varphi x_2^*x_2=x_2^*ux_1\varphi x_1^*u^*x_2
\;\;\;.
\]
As the left hand side is non-zero both sides are non-zero.
But $x_2^*x_2\varphi x_2^*x_2$ defines a representation equivalent
to $\rho_0$ and $x_2^*ux_1\varphi x_1^*u^*x_2$ defines a
representation equivalent to $\rho_{\gamma_1-\gamma_2}$.
But as these representations are disjoint the foregoing equality
gives a contradiction.
Thus all the non-zero representations $\rho_\gamma$ are disjoint.
Thus Condition~(4) implies Condition~(6) and the proposition is
established.\hfill$\Box$

\bigskip

\noindent{\bf Proof of Theorem \ref{tinvst3.1}}$\;$
It follows from the
introductory remarks to this section that we need to prove that if
$\omega$ is a pure state on $\gotA^\alpha$ , $x_1,x_2,\ldots,
x_n\in\gotA^\alpha$ and $\varepsilon>0$ then there is a pure state
$\varphi$ on $\gotA^\alpha$ satisfying any of the equivalent
conditions of Proposition \ref{pinvst3.1} such that
\[
|\varphi(x_i)-\omega(x_i)|<\varepsilon
\]
for $i=1,2,\ldots,n$.
So let $\omega, x_1,\ldots,x_n$ and $\varepsilon$ be given.
By Lemma \ref{linvst2.1} there is a positive element
$e\in\gotA^\alpha_{\sgotP}$ with $\omega(e)=1$ such that
\[
\|e\,x_i\,e-\omega(x_i)\,e^2\|<\varepsilon
\]
for $i=1,2,\ldots,n$.

Let $\{u_n\}$ be a dense sequence in the unitary group
$\gotU(\gotA^\alpha)$ of $\gotA^\alpha$, or $\gotA^\alpha+I$ if
$\gotA$ is non-unital.
Now $\gotA^\alpha(\gamma)\neq0$ for all
$\gamma\in\hat G$ since $G$ acts faithfully.
Let $\{\gamma_n\}$ be a sequence in
${\hat G}\backslash\{{\hat e}\}$ such that each element in
${\hat G}\backslash\{{\hat e}\}$ occurs infinitely often and
$\{(\gamma_n,u_n);\,n=1,2,\ldots\,\}$ is dense in
${\hat G}\backslash\{{\hat e}\}\times\gotU(\gotA^\alpha)$.
Next we construct inductively a sequence $e_n$ of
elements in $\gotA^\alpha_{\sgotP}$ and elements
$x_\gamma\in\gotA^\alpha(\gamma)$ with the three properties
\begin{description}
\item[{\rm(1)}]
$\;\;\;\;\;e_1\,e=e_1\;\;$ and
$\;\;e_n\,e_{n-1}=e_n\;\;\;\;n=2,3,\ldots$
\item[{\rm(2)}]
$\;\;\;\;\;e_n\,x_{\gamma_n}\,x_{\gamma_n}^*=e_n$
\item[{\rm(3)}]
$\;\;\;\;\;\|u_n\,x_{\gamma_n}^*\,e_n\,x_{\gamma_n}\,u_n^*\,e_n\|\leq1/n\;\;
\;.$
\end{description}
If these objects have been constructed for $1,2,\ldots,n-1$ we
construct $e_n$ and $x_{\gamma_n}$ as follows.
If $\gamma_n\in\{\gamma_1,\ldots,\gamma_{n-1}\}$ then
$x_{\gamma_n}$ has already been chosen and we keep that choice and
define $q_n=e_{n-1}$.
If $\gamma_n\not\in\{\gamma_1,\ldots,\gamma_{n-1}\}$ then we take
some  $q_n\in\gotA^\alpha_{\sgotP}$ with $q_n\,e_{n-1}=q_n$
and in addition some $x_{\gamma_n}\in\gotA^\alpha(\gamma_n)$ with
$q_n\,x_{\gamma_n}x_{\gamma_n}^*=q_n$.
This is possible since
$\gotA^\alpha(\gamma_n)\gotA^\alpha(\gamma_n)^*$ is an ideal in
$\gotA^\alpha$ which is prime by assumption.
Having now chosen
$x_{\gamma_n}, q_n$ choose $e_n$  in
$\gotA^\alpha_{\sgotP}$ with $e_n\,q_n=e_n$  such that
\[
\|u_n\,x_{\gamma_n}^*\,e_n\,x_{\gamma_n}\,u_n^*\,e_n\|\leq 1/n
\;\;\;\;.
\]
This choice is possible by the following reasoning.
Since $\gotA$ and  $\gotA^\alpha$ are prime and separable and
$\alpha_g$ is properly outer for all $g\in G\backslash\{e\}$ it
follows from \cite{BEEK} that there exists an irreducible
representation $\pi$ of $\gotA$ on a Hilbert space $\ch$ such that
$\pi(\gotA^\alpha)''=\gotB(\ch)$.
If {\it ad absurdum} $e_n$ does not exist then there is an
$\varepsilon>0$ such that
\[
\|u_n\,x_{\gamma_n}^*\,p\,x_{\gamma_n}\,u_n^*\,p\|\geq\varepsilon
\]
for all  $p\in\gotA^\alpha_{\sgotP}$ with $p\,q_n=p$.
But if $q$ is the spectral projection of $\pi(q_n)$ corresponding to
the eigenvalue one this means that
$|(\pi(xu_n^*)\psi,\psi)|\geq\varepsilon$ for all $\psi\in q\ch$
with $\|\psi\|=1$ where $x=x_{\gamma_n}$.
This conclusion is reached by noting that the orthogonal projection
$P_\psi$ onto $\psi$  is the limit of a decreasing net $\pi(p_m)$
with $p_m\,q_n=p_m$, $p_m\in\gotA^\alpha$, by irreducibility of $\pi$
on $\gotA^\alpha$.
But then one concludes that there is a $\theta\in\Ri$ such that
\[
\RRe(e^{i\theta}\pi(xu_n^*)\psi,\psi)\geq\varepsilon
\]
for all $\psi\in q\ch$ with $\|\psi\|=1$.
Thus
\[
e^{i\theta}fxu_n^*f+e^{-i\theta}fu_nx^*f\geq2\varepsilon f^2
\]
for any positive $f\in\gotA^\alpha$ with $fq_n=f$.
But applying $\alpha_g$ to this relation one obtains
\[
\gamma_n(g)e^{i\theta}fx_{\gamma_n}u_n^*f
+{\overline{\gamma_n(g)}}e^{-i\theta}fu_nx_{\gamma_n}^*f
\geq2\varepsilon f^2
\]
for all $g\in G$ since $u_n, f\in\gotA^\alpha$ and
$x_{\gamma_n}\in\gotA^\alpha(\gamma_n)$.
But this is impossible unless $\gamma_n=\hat e$.
Hence the $e_n$ exist.

\smallskip

Now one can complete the proof of Theorem \ref{tinvst3.1}.
Choose $e_n$ and $x_\gamma$ such that (1)--(3) are valid and
let $\varphi$ be any pure state on $\gotA^\alpha$ with
$\varphi(e_n)=1$ for all $n$.
Then $x_\gamma^*\varphi x_\gamma\disjoint \varphi$ for all
$\gamma\in{\hat G}\backslash\{\hat e\}$ by the following reasoning.
Let $p\in(\gotA^\alpha)^{**}$ be the closed projection such that
$e_n\searrow p$.
By Condition~(3) it follows that for any unitary $u\in\gotA^\alpha$
and any $\gamma\in{\hat G}\backslash\{\hat e\}$ one has
\[
u\,x_\gamma^*\,p\,x_\gamma \,u^*\,p=0
\]
but as $\varphi(p)=1$ this means that $x^*_\gamma\,\varphi\,
x_\gamma$ is disjoint from $\varphi$ for any
$\gamma\in{\hat G}\backslash\{\hat e\}$.
Now $\varphi\circ E$ is a pure state by Proposition \ref{pinvst3.1}.
But as
\[
\|e\,x_i\,e-\omega(x_i)\,e^2\|<\varepsilon
\]
and $p\,e=p$ one has
\[
\|p\,x_i\,p-\omega(x_i)\,p\|<\varepsilon
\]
and then
\[
|\varphi(x_i)-\omega(x_i)|=
|\varphi(p\,x_i\,p)-\omega(x_i)\varphi(p)|<\varepsilon
\;\;\;.
\]
This completes the proof of Theorem \ref{tinvst3.1} in the case of
abelian $G$.

Next we describe the necessary modification of the foregoing
argument to handle non-abelian $G$.
We first replace the dynamical system $(\gotA, G, \alpha)$ by
$({\widetilde{\gotA}}, G, {\tilde{\alpha}})$ where
\[
{\widetilde{\gotA}}=\gotA\otimes{\gotK}(L_2(G))\otimes \gotK\;\;\;,
\]
where ${\gotK}(L_2(G))$ denotes the compact operators on $L_2(G)$ and
$\gotK$ the compact operators on a separable infinite-dimensional
Hilbert space, and
\[
{\tilde\alpha}_g=\alpha_g\otimes\Ad\lambda(g)\otimes\iota
\]
where $\lambda$ is the left regular representation of $G$ on
$L_2(G)$.
If $\gamma\in\hat G$ we identify $\gamma$ with one of its concrete
representations and define $\gotA^\alpha_1(\gamma)$ (following
\cite{Was}, \cite{BrE} ) as the set of row-matrices
$x=(x_1,\ldots,x_d)$, where $x_i\in\gotA$ with $d=d(\gamma)=$
dimension of $\gamma$, such that
\[
\alpha_g (x)=x\gamma(g)
\;\;\;
\]
The reason for the replacement of $(\gotA,\alpha)$ by
$({\widetilde{\gotA}},\tilde{\alpha})$ is the existence of an
$x(\gamma)\in M({\widetilde{\gotA}}^{\tilde{\alpha}}_1(\gamma))$ with
the properties $x(\gamma)^*x(\gamma)=1\otimes1_{d(\gamma)}$ and
$x(\gamma)x(\gamma)^*=1$ \cite{BEK}.
Thus the problem of proving Theorem \ref{tinvst3.1} for $(\gotA,
G,\alpha)$ reduces to proving the theorem for
$({\widetilde{\gotA}}, G, \tilde{\alpha})$ by the following
reasoning.
 Let $\omega$ be a given $\alpha$-invariant state on
$\gotA$ and  assume $x_1,\ldots,x_n\in\gotA$ and $\varepsilon>0$ are
specified.
Then let $e$ be an $\Ad\lambda(g)\otimes\iota$-invariant
one-dimensional projection in ${\gotK}(L_2(G))\otimes\gotK$ and
replace $x_1,\ldots,x_n$ by  $\tilde x_1=x_1\otimes e,\ldots, \tilde
x_n=x_1\otimes e$ and $1\otimes e$ in ${\widetilde{\gotA}}$ and
$\omega$ by $\tilde\omega=\omega\otimes\omega'$ where $\omega'$ is
the pure state on ${\gotK}(L_2(G))\otimes\gotK$ such that
$\omega'(e)=1$. If Theorem \ref{tinvst3.1} is valid for
$({\widetilde{\gotA}}, G, \tilde{\alpha})$ let $\tilde{\varphi}$ be a
pure $\tilde{\alpha}$-invariant state on ${\widetilde{\gotA}}$ such
that
\begin{eqnarray*}
|\tilde{\varphi}(\tilde{x}_i)-\tilde{\omega}(\tilde{x}_i)|
&<&\varepsilon/2\\[5pt]
|\tilde{\varphi}(I\otimes e)-\tilde{\omega}(I\otimes e)|
&<&\varepsilon/2\;\;\;.
\end{eqnarray*}
But then $\tilde{\varphi}(I\otimes e)>1-\varepsilon/2$ so the state
 $\varphi$ defined on $\gotA$ by
\[
\varphi(x)=\tilde\varphi(x\otimes e)/\tilde{\varphi}(I\otimes e)
\]
is an $\alpha$-invariant pure state and
\[
|\varphi(x_i)-\omega(x_i)|<\varepsilon
\]
for $i=1,\ldots,n$.
Thus to prove Theorem \ref{tinvst3.1} we may assume there exist
$x(\gamma)\in\gotA^\alpha_1(\gamma)$ with
$x(\gamma)^*x(\gamma)=I\otimes I_{d(\gamma)}$,
$x(\gamma)x(\gamma)^*=I$.
One now defines dual endomorphisms ${\hat\alpha}_\gamma\colon
\gotA^\alpha\mapsto\gotA^\alpha$ by
\[
{\hat\alpha}_\gamma(a)=
\sum^{d(\gamma)}_{i=1}x_i(\gamma)ax_i(\gamma)^*=
x(\gamma)(a\otimes I_d)x^*(\gamma)
\]
for $a\in\gotA^\alpha$.
These are indeed endomorphisms of $\gotA^\alpha$ and if $\varphi$
is a state on $\gotA^\alpha$ one establishes, as in Proposition
\ref{pinvst3.1}, that $\varphi\circ E$ is a pure state on $\gotA$
if and only if $\varphi$ is pure on $\gotA^\alpha$ and
$\varphi\circ{\hat\alpha}_\gamma$ is disjoint to $\varphi$ for all
$\gamma\in{\hat G}\backslash\{e\}$.
To construct states with this property one uses the assumption that
$\alpha_g$ is properly outer and $\gotA$ and $\gotA^\alpha$ are
prime, which implies ${\hat\alpha}_\gamma$ is properly outer
\cite{BEK}, and then one uses the earlier argument to find minimal
projections $p\in(\gotA^\alpha)''$ such that
\[
\|u{\hat\alpha}_\gamma(p)u^*p\|=0
\]
for all unitaries $u\in\gotA^\alpha+\Ci I$.\hfill$\Box$

\section{One-dimensional shifts}

Let $\gotA$ be a unital $C^*$-algebra.
Then $\gotA$ is nuclear if and only if there exist an
increasing sequence $\{k_n\}$ of positive integers and unital
completely positive maps $\sigma_n$ of $\gotA$ into $M_{k_n}$
and $\tau_n$ of $M_{k_n}$ into $\gotA$ such that
\[
\lim_{n\to\infty}\|\tau_n\circ\sigma_n(x)-x\|=0
\]
for all $a\in \gotA$.
We now derive a more specific characterization of prime
AF-algebras, which are obviously nuclear.
\begin{lemma}\label{linvst4.1}
Let $\gotA$ be a unital $C^*$-algebra.
The following conditions are equivalent:
\begin{description}
\item[{\rm(1)}]
$\gotA$ is a prime {\rm AF}-algebra.
\item[{\rm(2)}]
There exist a sequence $\{k_n\}$ of positive
integers, unital completely positive maps $\sigma_n$ of
$\gotA$ onto $M_{k_n}$,  unital completely positive maps
$\tau_n$ of $M_{k_n}$ into $\gotA$ such that
$\tau_n(M_{k_n})$ is a $C^*$-subalgebra of $\gotA$,
homomorphisms $\iota_n$ {\rm(}not necessarily unital
$\,${\rm)} of $M_{k_n}$ into $\gotA$ such that the restriction
of $\iota_n\circ\sigma_n$ to $\gotA\bigcap\iota_n(M_{k_n})'$
has the form $\iota_n(I_{k_n})\,\omega_n(\cdot)$ where
$\omega_n$ is a pure state on $\gotA\bigcap\iota_n(M_{k_n})'$,
$\sigma_n\circ\iota_n$ is the identity
map on $M_{k_n}$ and
\[
\lim_{n\to\infty}\|\tau_n\circ\sigma_n(x)-x\|=0 \]
for all $x\in \gotA$.
\end{description}
\end{lemma}
\proof\ (2)$\Rightarrow$(1)
Since $\tau_n(M_{k_n})$ is a $C^*$-subalgebra of $\gotA$,
Condition~(2) implies that $\gotA$ is an AF-algebra by the
local characterization of AF-algebras in \cite{Bra1}.
That $\gotA$ is also prime follows the consequence of
Condition~(2) given in Propostion~\ref{pinvst4.1} below
(see Remark~\ref{rinvst4.1}).

Let $\{A_n\}$ be an increasing sequence of
finite-dimensional $C^*$-algebras such that
$ \gotA\cong\lim A_n$.
Let $A_n=\sum^{k_n}_{i=1} A_{ni}$ where the
$ A_{ni}$ are full matrix algebras, and let $e_{ni}$ be
the identity of $ A_{ni}$.
Since $ \gotA$ is prime we may suppose that the reduction
$ A_n\ni x\mapsto xe_{n+1,1}$ is faithful for any $n$,
\cite{Bra1}.

Suppose that $ A_{n+1,1}$ is isomorphic to $M_{k_n}$ and
define $\iota_n$ to be an isomorphism of $M_{k_n}$ onto
$ A_{n+1,1}$.
Let $\omega_n$ be a pure state of $ \gotA\bigcap A_{n+1,1}'$
such that $\omega_n(e_{n+1,1})=1$ and define $\sigma_n$ by
\[
\sigma_n(x)=(id\otimes \omega_n)(e_{n+1,1}x e_{n+1,1})
\]
where $e_{n+1,1} A e_{n+1,1}$ is identified with
$M_{k_n}\otimes
\Big(e_{n+1,1} \gotA e_{n+1,1}\bigcap A_{n+1,1}'\Big)$
by using $\iota_n\otimes id$.
Thus $\sigma_n$ is a unital completely positive map of
$ \gotA$ onto $M_{k_n}$ with $\sigma_n\circ\iota_n=id$ and
$\iota_n\circ\sigma_n$ restricted to
$ \gotA\bigcap\iota_n(M_{k_n})'$ is $\iota_n(I_{k_n})$
multiplied by the pure state $\omega_n$.
Choose a projection
$\varphi_n$ of $ A_{n+1,1}$ onto $ A_ne_{n+1,1}$ and define
$\tau_n$ by $\varphi_n\circ\iota_n$ where $ A_ne_{n+1,1}$ is
identified with $A_n$, i.e.,
$\tau_n=\psi_n\circ\varphi_n\circ\iota_n$ where
$\psi_n\colon A_ne_{n+1,1}\mapsto A_n$ is the $^*$-isomorphism
given by $\psi_n(xe_{n+1,1})=x$ for $x\in A_n$.
Thus $\tau_n$ is a unital completely positive map and
satisfies $\tau_n\circ\sigma_n(x)=x$ for all $x\in A_n$.
\hfill$\Box$

\bigskip

\begin{prop}\label{pinvst4.1}
Let $A$ be a unital, prime {\rm AF}-algebra,
$\gotA=\otimes^\infty_{-\infty}A$ the infinite tensor product
of $A$ and $\alpha$ the shift automorphism of $\gotA$.
Then the $\alpha$-invariant pure states of $\gotA$ are dense
in the set of $\alpha$-invariant states of $\gotA$.
\end{prop}
\proof\
Let $k_n$, $\sigma_n$, $\tau_n$ and $\iota_n$ be chosen as in
Lemma \ref{linvst4.1} with $A$ in place of $\gotA$.
Then define completely positive maps ${\tilde{\sigma}}_n$ of
$\gotA=A^{\otimes\infty}$ into
$\gotM_n={M_{k_n}}^{\otimes\infty}$ and ${\tilde{\tau}}_n$ of
$\gotM_n$ into $\gotA$ by \[
{\tilde{\sigma}}_n={\sigma_n}^{\otimes\infty}\;\;\;\;\;,\;\;\;\;\;
{\tilde{\tau}}_n={\tau_n}^{\otimes\infty}\;\;\;.
\]
Note that ${\tilde{\sigma}}_n\circ\alpha=
\alpha\circ{\tilde{\sigma}}_n$ and
${\tilde{\tau}}_n\circ\alpha= \alpha\circ{\tilde{\tau}}_n$
where we use $\alpha$ to denote the shift
automorphism for both $\gotA$ and $\gotM_n$.

Let $e_n=\iota_n(I_{k_n})$ and regard the infinite tensor product
 ${\tilde e}_n$ of copies of $e_n$ as a closed projection of
the second dual $\gotA^{**}$. Then define a homomorphism
${\tilde\iota}_n$ of $\gotM_n$ into ${\tilde e}_n\gotA{\tilde
e}_n$ by ${\iota_n}^{\otimes\infty}$. It again follows that
${\tilde\iota}_n\circ\alpha=\alpha\circ{\tilde\iota}_n$.

Next let $\omega$ be an $\alpha$-invariant state of $\gotA$.
Since the weak$^*$-limit of
$\omega\circ{\tilde\tau}_n\circ{\tilde{\sigma}}_n$ as
$n\to\infty$ is $\omega$ one may approximate
$\omega$ by
$\alpha$-invariant pure states by the following argument:
As $\omega\circ{\tilde\tau}_n$ is an $\alpha$-invariant
state of $\gotM_n$ one can find an $\alpha$-invariant pure
state $\varphi$ of $\gotM_n$ such that $\varphi$ is in a given
weak$^*$ neighbourhood of $\omega\circ{\tilde\tau}_n$
\cite{FNW}.
Then $\varphi\circ{\tilde{\sigma}}_n$ is in a given
neighbourhood of
$\omega\circ{\tilde\tau}_n\circ{\tilde{\sigma}}_n$.
Thus to complete the proof it suffices to show that
$\varphi\circ{\tilde{\sigma}}_n$ is pure.

Since
$\varphi\circ{\tilde{\sigma}}_n=
\varphi\circ{\tilde{\sigma}}_n\circ{\tilde{\iota}}_n
\circ{\tilde{\sigma}}_n$, the support of
$\varphi\circ{\tilde{\sigma}}_n$ is contained in
${\tilde e}_n$, and
\[
\varphi\circ{\tilde{\sigma}}_n
\bigg|_{(A\bigcap\iota_n(M_{k_n})')^{\otimes\infty}}
={\omega_n}^{\otimes\infty}
\;\;\;.
\]
Note also that
$\varphi\circ{\tilde{\sigma}}_n\circ{\tilde{\iota}}_n=\varphi$.
These facts determine $\varphi$ as follows:
for $x_i\in M_{k_n}$ and
$a_i\in e_n(A\bigcap\iota_n(M_{k_n})')e_n$ one has
\[
\varphi\circ{\tilde{\sigma}}_n
(\iota_n(x_1)a_1\otimes\iota_n(x_2)a_2\otimes\ldots
\iota_n(x_k)a_k)=
\varphi(x_1\otimes x_2\otimes\ldots\otimes x_k)
\omega_n(a_1)\omega_n(a_2)\ldots \omega_n(a_k)
\;\;\;,
\]
which is not ambiguous since $\varphi$ and
$\varphi\circ{\tilde{\sigma}}_n$ are $\alpha$-invariant.
Hence if $\psi$ is a positive linear functional on $\gotA$
such that $\psi\leq \varphi\circ{\tilde{\sigma}}_n$ then it
follows first that the support of $\psi$ is contained in
${\tilde e}_n$, and secondly, since $\omega_n$ is pure,
\[
\psi\bigg|_{(A\bigcap\iota_n(M_{k_n})')^{\otimes\infty}}
=\psi(I) \,{\omega_n}^{\otimes\infty}\;\;\;,
\]
and finally, since $\varphi$ is pure,
\[
\psi\circ{\tilde\iota}_n=\psi(I)\,\varphi\;\;\;.
\]
Hence it follows that
$\psi=\psi(I)\,\varphi\circ{\tilde{\sigma}}_n$.
Therefore $\varphi\circ{\tilde{\sigma}}_n$ is
pure.\hfill$\Box$

\medskip

\begin{remarkn}\label{rinvst4.1} If $A$ is not prime then the
conclusion of Proposition \ref{pinvst4.1} does not hold
because there are then two non-zero ideals $I$ and $J$ of $A$
such that $I\bigcap J=\{0\}$.
If we denote by $\hat I$ (respectively $\hat J$) the ideal of
$\gotA$ generated by $I$ (respectively $J$) at the position
$0$ in the factorization $\gotA=\otimes^\infty_{-\infty}A$
there is an $\alpha$-invariant state $\omega$ of $\gotA$
such that $\omega|_{\hat I}\neq0$ and $\omega|_{\hat J}\neq0$.
But since ${\hat I}\bigcap{\hat J}=\{0\}$, it follows that
$\varphi|_{\hat I}=0$ or $\varphi|_{\hat J}=0$ for any pure
state $\varphi$ on $\gotA$.
Thus $\omega$ cannot be approximated by a pure invariant
state, and not even by a pure state.
\end{remarkn}

\noindent{\large{\bf Acknowledgements}}$\;$ This work was carried out
whilst the first and third named authors were guests of the Mathematics
Department of Hokkaido University.
Their visits were sponsored by the  Norwegian Research
Council, the Australian Academy of Science and the Japan Society for
the Promotion of Science.


\end{document}